\documentclass[12pt]{article}
\usepackage{amsmath,amsthm,amsfonts,amssymb,amscd}

\allowdisplaybreaks[4]
\newtheorem {Lemma}{Lemma}[section]
\newtheorem {Theorem} {Theorem}[section]
\newtheorem{Proposition}{Proposition}[section]

\newenvironment {Proof} {\noindent {\bf Proof.}}{\quad $\square$\par\vspace{3mm}}

\begin{document}

\title{Ordering trees having small reverse Wiener indices}

\author{Rundan Xing, Bo Zhou\footnote{Corresponding author.}
\\
Department of Mathematics, South China Normal University, \\
Guangzhou 510631, P. R. China\\
email: {\tt zhoubo@scnu.edu.cn}
}

\date{}
\maketitle

\begin{abstract}
The reverse Wiener index of a connected graph $G$ is a variation
of the well-known Wiener index $W(G)$ defined as the sum of
distances between all unordered pairs of vertices of $G$. It is defined as
$\Lambda(G)=\frac{1}{2}n(n-1)d-W(G)$,
where $n$ is the number of vertices, and $d$ is the diameter of $G$.
We now determine the second and the third smallest reverse Wiener
indices of $n$-vertex trees and characterize the trees whose reverse
Wiener indices attain these values for $n\ge 6$ (it has been known that
the star is the unique tree with the smallest reverse Wiener index).
\end{abstract}

\section{Introduction}

Let $G$ be simple connected graph.
The Wiener index $W(G)$ of  $G$ is defined as the sum of
distances between all unordered pairs of vertices of $G$
\cite{GP,Ho,Wie}. The Wiener index is one of the oldest and the most useful
molecular--graph--based structure--descriptors used to explain
various chemical and physical properties of molecules and to
correlate the structure of molecules with their biological activity
\cite{Kar,Rou,ToCo,Tr}. It was also independently studied because of
its applications in  social science, architecture, and graph theory
\cite{Pl}. See \cite{GKM1,GKM2} for more details and results.
Its mathematical properties for trees can be found in
the review \cite{DoEnG} and in the references cited therein.

In 2000, Balaban {\it et al.\/}~\cite{BaMiIB} proposed a variant of
the Wiener index,  the reverse Wiener index. The reverse Wiener index of $G$ is
defined as~\cite{BaMiIB}
\[
\Lambda(G)=\frac{1}{2}n(n-1)d-W(G),
\]
where $n$ is the number of vertices and $d$ is the diameter of $G$.
The reverse Wiener index is also a useful structure--descriptor,
with applications in QSPR investigations, as demonstrated in
\cite{BaMiIB,IIB}.
Zhang and Zhou \cite{ZhZh} showed that the path $P_n$ and the star
$S_n$ are respectively the unique $n$-vertex trees with the largest
and the smallest reverse Wiener indices.
Cai and Zhou \cite{CaZh} determined the trees with the largest
reverse Wiener index within some subclasses of trees. Luo and Zhou
determined in \cite{LZ1} the $n$-vertex trees for $n\ge 5$  with the
$k$-th largest reverse Wiener indices for all $k$ up to
$\lfloor\frac{n}{2}\rfloor+1$, and determined in \cite{LZ2} the
$n$-vertex trees that are not caterpillars (trees for which the
removal of pendant vertices result in a path) for $n\ge 8$  with the
first a few largest reverse Wiener indices. See the survey \cite{Zh} for more results on the reverse Wiener index.

In this paper, we determine the second and the third smallest
reverse Wiener indices of $n$-vertex trees and characterize the
trees whose reverse Wiener indices attain these values for $n\ge 6$.

\section{Preliminaries}

Let $T$ be a tree  with   edge set $E(T)$. For
any $e\in E(T)$, $n_{T, 1}(e)$ and $n_{T, 2}(e)$ denote the number
of vertices of $T$ lying on the two sides of the edge $e$. For a
long time it has been known~\cite{GP,Wie}  that
\[
W(T)=\sum_{e\in E(T)}n_{T, 1}(e)\cdot n_{T, 2}(e).
\]

A vertex is pendent if it is of degree one. Let $\mathcal{T}_n$ be
the set of the trees on $n$ vertices. Let $\mathcal{CT}_n$ be the
set of trees $T$ in $\mathcal{T}_n$ such that one center of $T$ has
at least one pendent neighbor. We need the following two lemmas
given in \cite{ZhZh}. To be more self-contained in this paper,  we
include their proofs.

\begin{Lemma}   \label{CPT}
Let $T \in \mathcal{CT}_n$ with diameter $d \ge 4$. Then there is a
tree $T^* \in \mathcal{T}_n \setminus \mathcal{CT}_n$ with the same
diameter as $T$ such that $\Lambda(T^*) < \Lambda(T)$.
\end{Lemma}

\begin{Proof} Let $v$ be a center of $T$ with at least one pendent neighbor, say $w$. Obviously, there is at least one subtree $T_{1}$ in $T-v$  not containing  $w$ possesses
$n_{1}\le\frac{n}{2}-1$ vertices. Let  $v_1$ be the neighbor of $v$ in $T_1$. Let $T'$ be the tree formed from $T$ by
deleting edge $vw$ and adding edge $wv_1$. Note that both $T$ and $T'$ have diameter $d$. Then
\[
\Lambda(T')-\Lambda(T)=W(T)-W(T')=n_{1}(n-n_{1})-(n_{1}+1)(n-n_{1}-1)<0,
\]
and thus  $\Lambda (T')<\Lambda (T)$. Iterating the
transformation from $T$ to $T'$ will finally yield the tree $T^{*}$ as
required.
\end{Proof}

\begin{Lemma}   \label{d-2}
Let $T \in \mathcal{T}_n \setminus \mathcal{CT}_n$ with diameter $d
\ge 4$. Then there is a tree $T^* \in \mathcal{CT}_n$ with diameter
$d-2$ such that $\Lambda(T^*) < \Lambda(T)$.
\end{Lemma}

\begin{Proof}
Let  $v$ be a center of $T$  with neighbors
$v_1,v_2, \dots, v_p$.  For $i=1,2,\dots,p$,
let $T_{i}$ be the subtree in $T-v$ with $v_i\in V(T_i)$,
 $v_{i1}, v_{i2}, \dots, v_{ir_i}$ be the neighbors of $v_i$
in $T_i$, and $n_{i}=|V(T_i)|$.
 Let  $T^{*}$ be the tree formed from
$T$ by deleting edges $v_iv_{ij}$ and adding edges $vv_{ij}$ for
all $i=1,2 \dots, p$ and $j=1, 2,\dots, r_i$. Obviously,  $T^* \in \mathcal{CT}_n$.
It is easily seen that
\[
\renewcommand{\arraystretch}{1.8}
\begin{array}{lll}
\Lambda (T^{*})-\Lambda (T)&=&-n(n-1)+W(T)-W(T^{*})\\ &=&
-n(n-1)+\sum\limits_{i=1}^{p}n_{i}(n-n_{i})-p(n-1)
\\
&\leq& -n(n-1)+(n-1)\sum\limits_{i=1}^{p}n_i-p(n-1)\\
&=& -n(n-1)+(n-1)(n-1)-p(n-1)\\
&=&-(n-1)(1+p)<0,
\end{array}
\]
as desired.
\end{Proof}

An $n$-vertex tree of diameter $4$ may be constructed as follows:
For some integer $k\ge 2$, it is obtained from the star $S_{k+1}$
with center $v_0$ and pendent vertices $v_1,v_2,\dots,v_k$ by
attaching $m_i$ pendent vertices to $v_i$ for $i=0,1,\ldots,k$,
where $m_0 \ge 0$, $m_i \ge 1$ for $i=1,\ldots,k$, and
$\sum_{i=0}^km_i +k+1=n$. If there are $s$ distinct (positive)
numbers $n_1 < n_2 < \cdots < n_s$ among the numbers
$m_1,m_2,\ldots,m_k$, where $n_i$ appears $b_i$ times for
$i=1,2,\ldots,s$, then such a tree is denoted by
$T_{n,k}\left(n_0;n_1^{[b_1]},n_2^{[b_2]},\ldots,n_s^{[b_s]}\right)$
with $n_0=m_0$, where $k \ge 2, \sum_{i=1}^{s}b_i=k$ and
$\sum_{i=1}^{s}{b_i n_i}=n-n_0-k-1$. If $n_0=0$, then we write it
simply
$T_{n,k}\left(n_1^{[b_1]},n_2^{[b_2]},\ldots,n_s^{[b_s]}\right)$.


Obviously, for an integer $n\ge 2$, there is a unique positive
integer $q$ such that $q^2<n\le (q+1)^2$, and thus $n$ may be
written as
\begin{equation} \label{formu-n}
n=q^2+r \mbox{ with }r=1,2,\ldots,2q+1.
\end{equation}
It is easily seen that
\begin{eqnarray*}
& &
   W\left(T_{n,k}\left(n_0;n_1^{[b_1]},n_2^{[b_2]},\ldots,n_s^{[b_s]}\right)\right)\\
&=&
   1\cdot(n-1)\cdot (n-1-k)+\sum_{i=1}^{s}{b_i(n_i+1)(n-n_i-1)}\\
&=&
   (n-1)^2-(n-1)k+n\sum_{i=1}^{s}{b_i(n_i+1)}-\sum_{i=1}^{s}{b_i(n_i+1)^2}\\
&=&
   (n-1)^2-(n-1)k+n(n-n_0-1)-2(n-n_0-k-1)-k-\sum_{i=1}^{s}{b_in_i^2}\\
&=&
   (n-1)(2n-3)-(n-2)k-(n-2)n_0-\sum_{i=1}^{s}{b_in_i^2}.
\end{eqnarray*}

\begin{Lemma} \label{balace}
Let
$T=T_{n,k}\left(n_0;n_1^{[b_1]},n_2^{[b_2]},\ldots,n_s^{[b_s]}\right)$
with $k \ge 2$. If  $T$ achieves the smallest reverse Wiener index
for fixed $n$, then $s \le 2$, and  $n_2-n_1 = 1$ if $s=2$.
\end{Lemma}

\begin{Proof}
Suppose that there exist $i$ and $j$ with $1 \le i,j \le s$ such
that $n_i-n_j \ge 2$ (which is obviously true if $s\ge 3$).  Suppose
without loss of generality that $v_{i'}$ has $n_i$ pendent
neighbors, one of which is denoted by $u$, and $v_{j'}$ has $n_j$
pendent neighbors, where $1\le i',j'\le k$. Let
$T^*=T-uv_{i'}+uv_{j'}$. Then
\begin{eqnarray*}
   \Lambda(T)-\Lambda(T^*)
&=&
   W(T^*)-W(T)\\
&=&
   n_i(n-n_i)+(n_j+2)(n-n_j-2)\\
&&-\left[(n_i+1)(n-n_i-1)+(n_j+1)(n-n_j-1)\right]\\
&=&
   2(n_i-n_j-1) > 0,
\end{eqnarray*}
and thus  $\Lambda(T)>\Lambda(T^*)$, a contradiction. The result
follows.
\end{Proof}

\begin{Lemma} \label{k^2+1}
For fixed $n$ of the form $(\ref{formu-n})$, if
$T_{n,k}\left(n_1^{[b_1]},n_2^{[b_2]},\ldots,n_s^{[b_s]}\right)$
with $k \ge 2$ and $s\le 2$ where $n_2-n_1=1$ if $s=2$ achieves the
smallest reverse Wiener index for fixed $n$, then $k=q-1$, $q$ if
$r=1$, and $k=q$, $q+1$ if $r=2,\ldots,2q+1$.
\end{Lemma}

\begin{Proof}
Let
$T=T_{n,k}\left(n_1^{[b_1]},n_2^{[b_2]},\ldots,n_s^{[b_s]}\right)$,
where $s=1,2$, and $n_2-n_1 = 1$ if $s=2$, be a tree with the
smallest reverse Wiener index for fixed $n$.

If $s=1$, then  $n_1=\frac{n-1}{k}-1$, $b_1=k$, and thus
\[
\Lambda(T)=2n(n-1)-W(T)=n-1+(n-1)k+\frac{(n-1)^2}{k}.
\]
If $s=2$, then $n_1=\left\lfloor\frac{n-1}{k}\right\rfloor-1$,
$n_2=\left\lfloor\frac{n-1}{k}\right\rfloor$,
$b_1=\left\lfloor\frac{n-1}{k}\right\rfloor k-(n-k-1)$,
$b_2=n-1-\left\lfloor\frac{n-1}{k}\right\rfloor k$, and thus
\begin{eqnarray*}
   \Lambda(T)&=&2n(n-1)-W(T)\\
&=&
   3(n-1)+(n-2)k+\left(\left\lfloor\frac{n-1}{k}\right\rfloor
   k-n+k+1\right)\cdot\left(\left\lfloor\frac{n-1}{k}\right\rfloor-1\right)^2\\
& &
   +\left(n-1-\left\lfloor\frac{n-1}{k}\right\rfloor
   k\right)\cdot\left\lfloor\frac{n-1}{k}\right\rfloor^2.
\end{eqnarray*}
Suppose that $n-1 \equiv t \pmod k$, where $t=0,1,2,\ldots,k-1$.
Then $\left\lfloor\frac{n-1}{k}\right\rfloor=\frac{n-t-1}{k}$, and
thus
\[
\Lambda(T)= n+t-1+(n-1)k+\frac{(n-1)^2-t^2}{k}.
\]
Let $f(k)$ be the expression in the right side of the equation
above. Then
\[
f'(k)=(n-1)-\frac{(n-1)^2-t^2}{k^2},
\]
from which we know that $f(k)$ is decreasing for $k \le
\sqrt{(n-1)-\frac{t^2}{n-1}}$, and increasing for $k \ge
\sqrt{(n-1)-\frac{t^2}{n-1}}$. If $r=1$, then on one hand,
$\sqrt{(n-1)-\frac{t^2}{n-1}}=\sqrt{q^2-\frac{t^2}{q^2}}\le q$,
implying that $k\le q$, and on the other hand,
$\sqrt{q^2-\frac{t^2}{q^2}}\ge \sqrt{q^2-\frac{(k-1)^2}{q^2}}$,
implying that $k\ge \sqrt{q^2-\frac{(k-1)^2}{q^2}}$ and then $k\ge
q-1$. If $r\ge 2$, then on one hand, $\sqrt{n-1-\frac{t^2}{n-1}} \le
\sqrt{n-1} \le \sqrt{q^2+2q} < q+1$, implying that $k\le q+1$, and
on the other hand, $\sqrt{n-1-\frac{t^2}{n-1}}\ge
\sqrt{n-1-\frac{q^2}{n-1}}>q$, implying that $k\ge q$. The result
follows.
\end{Proof}

Let $\mathbb{T}_{n,d}$ be the set of $n$-vertex trees with diameter
$d$, where $2 \le d \le n-1$.

\begin{Lemma} \label{T-5}
Let $T$ be a tree in $\mathbb{T}_{n,5}$ with the smallest reverse
Wiener index. Then there exists $T^*\in \mathcal{CT}_n$ with
diameter $4$ such that $\Lambda(T^*) < \Lambda(T)$.
\end{Lemma}

\begin{Proof}
Let $T$ be a tree in $\mathbb{T}_{n,5}$ with the smallest reverse
Wiener index. By Lemma \ref{CPT}, there exists no pendent neighbor
for the centers of $T$. Let $T^*$ be the tree obtained by
contracting the edge connecting the centers  $u$ and $v$ of $T$
followed by attaching a pendent vertex to the new vertex  resulting
from identifying $u$ and $v$. Then $T^* \in \mathbb{T}_{n,4}$ and
its center has a pendent neighbor. It is easily seen that
\begin{eqnarray*}
& &
   \Lambda(T)-\Lambda(T^*)\\
&=&
   \frac{5}{2}n(n-1)-W(T)-\left(2n(n-1)-W(T^*)\right)\\
&=&
   \frac{5}{2}n(n-1)-n_{T,1}(uv) n_{T,2}(uv)-2n(n-1)+(n-1)\\
&\ge&
   \frac{1}{2}(n-1)(n+2)-\left\lfloor\frac{n}{2}\right\rfloor
   \left\lceil\frac{n}{2}\right\rceil >0,
\end{eqnarray*}
from which we get the desired result.
\end{Proof}

For $2 \le d \le n-1$, let $f(n,d)=\min \{\Lambda(T): T\in
\mathbb{T}_{n,d}\}$, and let $\mathcal{T}_{n,d}$ be the set of trees
in $\mathbb{T}_{n,d}$ with reverse Wiener index $f(n,d)$.

Obviously, $\mathbb{T}_{n,2}=\{S_n\}$, and
$f(n,2)=\Lambda(S_n)=n-1$.

For $3 \le d \le n-2$, let $g(n,d)=\min \{\Lambda(T): T\in
\mathbb{T}_{n,d}\setminus \mathcal{T}_{n,d}\}$, and
$\mathcal{T'}_{n,d}$ be the set of trees in $\mathbb{T}_{n,d}$ with
reverse Wiener index $g(n,d)$.

 For  $2 \le a \le
\lfloor\frac{n}{2}\rfloor$, let $D_{n,a}$ be the tree formed by
adding an edge between the centers of the stars $S_a$ and $S_{n-a}$.

\begin{Proposition} \label{d=3}
For $n\ge 4$,
\[
f(n,3)=\frac{n^2}{2}+\frac{3n}{2}-2-\left\lfloor\frac{n}{2}\right\rfloor
\left\lceil\frac{n}{2}\right\rceil \mbox{ and }
\mathcal{T}_{n,3}=\{D_{n,\lfloor\frac{n}{2}\rfloor}\}.
\]
For $n\ge 6$,
\[
g(n,3)=\frac{n^2}{2}+\frac{3n}{2}-2-\left(\left\lfloor\frac{n}{2}\right\rfloor-1\right)
\left(\left\lceil\frac{n}{2}\right\rceil+1\right) \mbox{ and
}\mathcal{T'}_{n,3}=\{D_{n,\lfloor\frac{n}{2}\rfloor-1}\}.
\]
%
\end{Proposition}

\begin{Proof}
Obviously, the trees in $\mathbb{T}_{n,3}$ are of the type
$D_{n,a}$, where $2 \le a \le \lfloor\frac{n}{2}\rfloor$. For $2 \le
a \le \lfloor\frac{n}{2}\rfloor$, we have
\begin{eqnarray*}
\Lambda(D_{n,a}) &=&
   \frac{3}{2}n(n-1)-W(D_{n,a})\\
&=&
   \frac{3}{2}n(n-1)-\left[(n-1)(n-2)+a(n-a)\right]\\
&=&
   \frac{n^2}{2}+\frac{3n}{2}-2-a(n-a),
\end{eqnarray*}
from which we know that, for fixed $n$, $\Lambda(D_{n,a})$ is
decreasing for $2 \le a \le \lfloor\frac{n}{2}\rfloor$. The result
follows.
\end{Proof}

\section{Trees with the second smallest reverse Wiener index}

In this section, we determine the second smallest reverse Wiener
index of $n$-vertex trees and characterize the trees achieving this
value.

\begin{Proposition} \label{T_4} Let $n\ge 5$ be of the form
$(\ref{formu-n})$. Then
\begin{eqnarray*}
f(n,4)=
\begin{cases}
2q^3+q^2+3rq-3q+2r-2&    \ \     \mbox{if $r=1,\ldots,q$,}\\
2q^3+q^2+3rq-4q+3r-3&    \ \     \mbox{if $r=q+1,\dots,2q+1$,}
\end{cases}
\end{eqnarray*}

\begin{eqnarray*}
\mathcal{T}_{n,4}=
\begin{cases}
\left\{T_{n,q}\left(q-1^{[q]}\right)\right\}&    \ \     \mbox{if $r=1$,}\\
\left\{T_{n,q}\left(q-1^{[q-r+1]},q^{[r-1]}\right)\right\}&    \ \     \mbox{if $r=2,\ldots,q$,}\\
\left\{T_{n,q}\left(q^{[q]}\right),T_{n,q+1}\left(q-1^{[q+1]}\right)\right\}&    \ \     \mbox{if $r=q+1$,}\\
\left\{T_{n,q+1}\left(q-1^{[2(q+1)-r]},q^{[r-q-1]}\right)\right\}& \
\ \mbox{if $r=q+2,\ldots,2q+1$.}
\end{cases}
\end{eqnarray*}
\end{Proposition}

\begin{Proof}
Let $T$ be a tree with the smallest reverse Wiener index in
$\mathbb{T}_{n,4}$. By Lemma \ref{CPT}, $T$ must be written as
$T=T_{n,k}\left(n_1^{[b_1]},n_2^{[b_2]},\ldots,n_s^{[b_s]}\right)$
with $k \ge 2$ and $\sum_{i=1}^{s}{b_i n_i}=n-k-1$. By Lemmas
\ref{balace} and \ref{k^2+1}, we have $s \le 2$, $k=q-1$, $q$ if
$r=1$, and $k=q$, $q+1$ if $r\ge 2$.

Suppose that $k=q$. If $s=1$, then
$n=q^2+r=1+k+n_1k=1+q+n_1q$, and thus $r-1=q(1+n_1-q)$, which implies
that $r=1,q+1,2q+1$. If $s=2$, then
$n=q^2+r=1+q+b_1n_1+b_2n_2=1+q+b_1n_1+b_2(n_1+1)=1+q+qn_1+b_2$ where
$1 \le b_2 \le q-1$, and thus $r-b_2-1=q(n_1+1-q)$, which implies
that $r=2,\ldots,q,q+2,\ldots,2q$, i.e., $r \ne 1,q+1,2q+1$. In
conclusion, if $k=q$, then $s=1$ if and only if
$r=1,q+1,2q+1$.
Now suppose that $k=q+1$. If $s=1$, then
$n=q^2+r=1+(q+1)+n_1(q+1)$, and thus $r=(q+1)(n_1+2-q)$, which
implies that $r=q+1$. If $s=2$, then
$n=q^2+r=1+(q+1)+b_1n_1+b_2n_2=q+2+(q+1)n_1+b_2$ where $1 \le b_2
\le q$, and thus $r-b_2=(q+1)(n_1+2-q)$, which implies that
$r=1,\ldots,q,q+2,\ldots,2q+1$, i.e., $r \ne q+1$. In conclusion, if
$k=q+1$, then $s=1$  if and only if $r=q+1$.

\noindent {\bf Case 1.} $r=1$. Note that $k=q-1$, $q$. If $k=q-1$,
then $s=2$, and thus $T=T_{n,q-1}\left(q^{[q-2]}, q+1^{[1]}\right)$
with $\Lambda(T)=2q^3+q^2+q+2$. If $k=q$, then  $s=1$, and
$T=T_{n,q}\left(q-1^{[q]}\right)$ with
$\Lambda(T)=2q^3+q^2<2q^3+q^2+q+2$. Thus
$T=T_{n,q}\left(q-1^{[q]}\right)$.

\noindent {\bf Case 2.} $r=2,\ldots,q$. If $k=q$, then $s=2$, and
thus $n_1=\left\lfloor\frac{q^2+r-1}{q}\right\rfloor-1=q-1$,
$n_2=q$, $b_1=q^2-(q^2+r-q-1)=q-r+1$ and $b_2=r-1$. If  $k=q+1$,
then $s=2$, and thus $n_1=q-2$, $n_2=q-1$, $b_1=q-r+1$ and $b_2=r$.
Thus $T=T_{n,q}\left(q-1^{[q-r+1]},q^{[r-1]}\right)$,
$T_{n,q+1}\left(q-2^{[q-r+1]},q-1^{[r]}\right)$. By direct
calculation,
\begin{eqnarray*}
&&
   \Lambda\left(T_{n,q}\left(q-1^{[q-r+1]},q^{[r-1]}\right)\right)\\
&=&
   2n(n-1)-[(n-1)(n-1-q)+q(n-q)(q-r+1)\\
& &
   (q+1)(n-q-1)+(r-1)]\\
&=&
   n^2+n-nq^2+nq-nr+q^3+2rq-3q+r-2\\
&=&
   2q^3+q^2+3rq-3q+2r-2,
\end{eqnarray*}
\begin{eqnarray*}
&&
   \Lambda\left(T_{n,q+1}\left(q-2^{[q-r+1]},q-1^{[r]}\right)\right)\\
&=&
   2n(n-1)-[(n-1)(n-2-q)\\
& &
   +(q-1)(n-q+1)(q-r+1)+q(n-q)r]\\
&=&
   n^2+2n-nq^2+nq-nr+q^3-q^2+2rq-2q-r-1\\
&=&
   2q^3+q^2+3rq-2q+r-1.
\end{eqnarray*}
Since
$(2q^3+q^2+3rq-3q+2r-2)-(2q^3+q^2+3rq-2q+r-1)=
r-q-1 \le -1 <0$, we have
 $\Lambda\left(T_{n,q}\left(q-1^{[q-r+1]},q^{[r-1]}\right)\right)
<\Lambda\left(T_{n,q+1}\left(q-2^{[q-r+1]},q-1^{[r]}\right)\right)$, and then
$T=T_{n,q}\left(q-1^{[q-r+1]},q^{[r-1]}\right)$.

\noindent {\bf Case 3.} $r=q+1$. By direct calculation,
$\Lambda\left(T_{n,q}\left(q^{[q]}\right)\right) =
\Lambda\left(T_{n,q+1}\left(q-1^{[q+1]}\right)\right)
=2q^3+4q^2+2q$. Since $s=1$, we have
$T=T_{n,q}\left(q^{[q]}\right)$,
$T_{n,q+1}\left(q-1^{[q+1]}\right)$.

\noindent {\bf Case 4.} $r=q+2,\ldots,2q$. If $k=q$, then $s=2$, and
thus $n_1=q$, $n_2=q+1$, $b_1=2q+1-r$ and $b_2=r-q-1$. If $k=q+1$,
then $s=2$, and thus $n_1=q-1$, $n_2=q$,  $b_1=2(q+1)-r$ and
$b_2=r-q-1$. Thus
$T=T_{n,q}\left(q^{[2q+1-r]},q+1^{[r-q-1]}\right)$,
$T_{n,q+1}\left(q-1^{[2(q+1)-r]},q^{[r-q-1]}\right)$. By direct
calculation,
\begin{eqnarray*}
&&
   \Lambda\left(T_{n,q}\left(q^{[2q+1-r]},q+1^{[r-q-1]}\right)\right)\\
&=&
   2n(n-1)-[(n-1)(n-1-q)+(q+1)(n-q-1)(2q+1-r)\\
& &
   +(q+2)(n-q-2)(r-q-1)]\\
&=&
   n^2+n-nq^2+nq-nr+q^3+2rq-5q+3r-4\\
&=&
   2q^3+q^2+3rq-5q+4r-4,
\end{eqnarray*}
\begin{eqnarray*}
&&
   \Lambda\left(T_{n,q+1}\left(q-1^{[2(q+1)-r]},q^{[r-q-1]}\right)\right)\\
&=&
   2n(n-1)-[(n-1)(n-2-q)+q(n-q)(2q+2-r)\\
   &&+(q+1)(n-q-1)(r-q-1)]\\
&=&
   n^2+2n-nq^2+nq-nr+q^3-q^2+2rq-4q+r-3\\
&=&
   2q^3+q^2+3rq-4q+3r-3.
\end{eqnarray*}
Note that
$\Lambda\left(T_{n,q}\left(q^{[2q+1-r]},q+1^{[r-q-1]}\right)\right)
-
\Lambda\left(T_{n,q+1}\left(q-1^{[2(q+1)-r]},q^{[r-q-1]}\right)\right)
= r-q-1 \ge 1>0$. Then $T=T_{n,q+1}\left(q-1^{[2(q+1)-r]},q^{[r-q-1]}\right)$.

\noindent {\bf Case 5.}  $r=2q+1$. If $k=q$, then $s=1$. If $k=q+1$,
then $s=2$. Thus $T=T_{n,q}\left(q+1^{[q]}\right)$,
$T_{n,q+1}\left(q-1^{[1]},q^{[q]}\right)$. Note that
$\Lambda\left(T_{n,q}\left(q+1^{[q]}\right)\right)=2q^3+7q^2+6q>
\Lambda\left(T_{n,q+1}\left(q-1^{[1]},q^{[q]}\right)\right)=2q^3+7q^2+5q$.
Then $T=T_{n,q+1}\left(q-1^{[1]},q^{[q]}\right)$.
\end{Proof}

Now we are ready to give our main result in this section.

\begin{Theorem} \label{sec-min}
Among the trees in $\mathcal{T}_n$ with $n \ge 4$,
$D_{n,\lfloor\frac{n}{2}\rfloor}$ for $n \le 56$,  $D_{57,28}$,
$T_{57,7}\left(7^{[7]}\right)$ and $T_{57,8}\left(6^{[8]}\right)$
for $n=57$, and the trees in $\mathcal{T}_{n,4}$ for $n \ge 58$ are
the unique trees with the second smallest reverse Wiener index,
which is equal to
$\frac{n^2}{2}+\frac{3n}{2}-2-\lfloor\frac{n}{2}\rfloor
\lceil\frac{n}{2}\rceil$ for $n \le 56$, $896$ for $n=57$, and
$f(n,4)$ for $n\ge 58$, respectively, where $\mathcal{T}_{n,4}$ and
$f(n,4)$ are given in Proposition \ref{T_4}.
\end{Theorem}

\begin{Proof} The case $n=4$ is trivial.  Suppose that $n\ge 5$.
Let $T\in \mathcal{T}_n$. Let $d$ be the diameter of $T$. If $d\ge 5$, then
by Lemmas \ref{CPT} and \ref{d-2}, we have
\begin{eqnarray*}
\Lambda(T)\ge f(n,d)>f(n,4)>f(n,2) && \mbox{for even } d,\\
\Lambda(T)\ge f(n,d) \ge f(n,5)>f(n,3) && \mbox{for odd }d.
\end{eqnarray*}
Thus the second smallest reverse Wiener index of the trees in
$\mathcal{T}_n$ is equal to $\min\{f(n,3), f(n,4)\}$, and by
Propositions \ref{d=3} and \ref{T_4},  it is only achieved by
$D_{n,\lfloor\frac{n}{2}\rfloor}$ or trees in $\mathcal{T}_{n,4}$,
where, with $n$ being of the form $(\ref{formu-n})$,
\begin{eqnarray*}
f(n,3)&=&
\begin{cases}
\frac{1}{4}q^4+\frac{1}{2}rq^2+\frac{3}{2}q^2+\frac{1}{4}r^2+\frac{3}{2}r-\frac{7}{4}&
\mbox{if $n$ is odd,}\\
\frac{1}{4}q^4+\frac{1}{2}rq^2+\frac{3}{2}q^2+\frac{1}{4}r^2+\frac{3}{2}r-2&
\mbox{if $n$ is even,}
\end{cases}
\end{eqnarray*}
and $f(n,4)$ is given in Proposition \ref{T_4}. If $n=5,6,\dots,56$,
then it can be checked that $f(n,3)<f(n,4)$, and thus the result
follows from Proposition \ref{d=3}.  Suppose that $n \ge 57$.

\noindent {\bf Case 1.} $r=1,2,\dots,q$. Since $q^2+1\ge 57$, we
have $q\ge 8$. Note that $f(n,4)=2q^3+q^2+3rq-3q+2r-2$. If $n$ is
odd, then since $\frac{1}{2}q^2-3q-\frac{1}{2} >0$, we have
\begin{eqnarray*}
& &
   f(n,3)-f(n,4)\\
&=&
   \frac{1}{4}q^4+\frac{1}{2}rq^2+\frac{3}{2}q^2+\frac{1}{4}r^2+\frac{3}{2}r-\frac{7}{4}
   -\left(2q^3+q^2+3rq-3q+2r-2\right)\\
&=&
   \frac{1}{4}q^4-2q^3+\frac{1}{2}q^2+3q+\frac{1}{4}+\frac{1}{4}r^2+\left(\frac{1}{2}q^2-3q-\frac{1}{2}\right)r\\
&\ge&
   \frac{1}{4}q^4-2q^3+\frac{1}{2}q^2+3q+\frac{1}{4}+\frac{1}{4}+\left(\frac{1}{2}q^2-3q-\frac{1}{2}\right)\\
&=&
   \frac{1}{4}q^4-2q^3+q^2>0,
\end{eqnarray*}
and if $n$ is even, then
\[
f(n,3)-f(n,4) \ge \frac{1}{4}q^4-2q^3+q^2-\frac{1}{4}>0.
\]
Thus $f(n,3)>f(n,4)$.

\noindent {\bf Case 2.} $r=q+1,q+2,\ldots,2q+1$. Since $q^2+q+1 \ge
57$, then $q\ge 7$. Note that $f(n,4)=2q^3+q^2+3rq-4q+3r-3$. Suppose
first that $n$ is odd. Since $\frac{1}{2}q^2-3q-\frac{3}{2}
>0$, we have
\begin{eqnarray*}
& &
   f(n,3)-f(n,4)\\
&=&
   \frac{1}{4}q^4+\frac{1}{2}rq^2+\frac{3}{2}q^2+\frac{1}{4}r^2+\frac{3}{2}r-\frac{7}{4}
   -\left(2q^3+q^2+3rq-4q+3r-3\right)\\
&=&
   \frac{1}{4}q^4-2q^3+\frac{1}{2}q^2+4q+\frac{5}{4}+\frac{1}{4}r^2+\left(\frac{1}{2}q^2-3q-\frac{3}{2}\right)r\\
&\ge&
   \frac{1}{4}q^4-2q^3+\frac{1}{2}q^2+4q+\frac{5}{4}+\frac{1}{4}(q+1)^2
   +\left(\frac{1}{2}q^2-3q-\frac{3}{2}\right)(q+1)\\
&=&
   \frac{1}{4}q^4-\frac{3}{2}q^3-\frac{7}{4}q^2 \ge 0,
\end{eqnarray*}
and then $f(n,3)\ge f(n,4)$ with equality if and only if $q=7$ and $r=q+1=8$, i.e., $n=57$. 
Now suppose that $n$ is even. If $q=7$, then $n=58,60,62,64$, and it
is easily checked by the expressions for $f(n,3)$ and $f(n,4)$ that
$f(n,3)>f(n,4)$. If $q\ge 8$, then
\[
f(n,3)-f(n,4) \ge
\frac{1}{4}q^4-\frac{3}{2}q^3-\frac{7}{4}q^2-\frac{1}{4} > 0,
\]
and thus $f(n,3)>f(n,4)$.

 Combining Cases 1 and 2,
$f(n,3)>f(n,4)$ for $n \ge 58$ and $f(57,3)=f(57,4)$.
The result for $n\ge 57$ follows from Propositions \ref{d=3} and
\ref{T_4}.
\end{Proof}

\section{Trees with the third smallest reverse Wiener index}

In this section, we determine the third smallest reverse Wiener
index of $n$-vertex trees and characterize the trees whose reverse
Wiener index achieve this value.

\begin{Proposition} \label{T_4-sec}  Let $n\ge 6$ be of the form
$(\ref{formu-n})$. Then
\begin{eqnarray*}
g(n,4)=
\begin{cases}
2q^3+q^2+3rq-3q+2r&    \ \     \mbox{if $r=1,\ldots,q-1$,}\\
2q^3+4q^2+5q+4&    \ \     \mbox{if $r=q+2$,}\\
2q^3+q^2+3rq-4q+3r-1&    \ \     \mbox{if
$r=q,q+1,q+3,\ldots,2q+1$,}
\end{cases}
\end{eqnarray*}
\[
\mathcal{T'}_{n,4} =\begin{cases}
\{T_{n,q}\left(q-2^{[1]},q-1^{[q-r-1]},q^{[r]}\right)\}
\mbox{~~~~~~}
\mbox{if $r=1,\ldots,q-2$,}\\

\{T_{n,q}\left(q-2^{[1]},q^{[q-1]}\right), T_{n,q+1}\left(q-2^{[2]},q-1^{[q-1]}\right)\} \\
\mbox{~~~~~~~~~~~~~~~~~~~~~~~~~~~~~~~~~~~~~~~~~~~~~}  \mbox{if $r=q-1$,}\\

\{T_{n,q+1}\left(q-2^{[1]},q-1^{[q]}\right)\}   \mbox{~~~~~~~~~~~~~~}  \mbox{if $r=q$,}\\

\{T_{n,q}\left(q-1^{[1]},q^{[q-2]},q+1^{[1]}\right),
T_{n,q+1}\left(q-2^{[1]},q-1^{[q-1]},q^{[1]}\right)\} \\   \mbox{~~~~~~~~~~~~~~~~~~~~~~~~~~~~~~~~~~~~~~~~~~~~~}
\mbox{if $r=q+1$,}\\

\{T_{n,q}\left(q^{[q-1]},q+1^{[1]}\right)\}   \mbox{~~~~~~~~~~~~~~~~~~~}  \mbox{if $r=q+2$,}\\

\{T_{n,q}\left(q^{[q-2]},q+1^{[2]}\right),T_{n,q+1}\left(q-2^{[1]},q-1^{[q-3]},q^{[3]}\right),\\
   T_{n,q+1}\left(q-1^{[q]},q+1^{[1]}\right)\}  \mbox{~~~~~~~~~~~~~~~~}  \mbox{if $r=q+3$,}\\

\{T_{n,q+1}\left(q-2^{[1]},q-1^{[2q-r]},q^{[r-q]}\right),\\
T_{n,q+1}\left(q-1^{[2q+3-r]},q^{[r-q-3]},q+1^{[1]}\right)\}\\
\mbox{~~~~~~~~~~~~~~~~~~~~~~~~~~~~~~~~~~~~~~~~~~~~~}
\mbox{if $r=q+4,\ldots,2q-1$,}\\

\{T_{n,q+1}\left(q-2^{[1]},q^{[q]}\right),
T_{n,q+1}\left(q-1^{[2q+3-r]},q^{[r-q-3]},q+1^{[1]}\right)\}\\
\mbox{~~~~~~~~~~~~~~~~~~~~~~~~~~~~~~~~~~~~~~~~~~~~~}
\mbox{if $r=2q$,}\\

\{T_{n,q+1}\left(q-1^{[2]},q^{[q-2]},q+1^{[1]}\right)\}
\mbox{~~~~~~}  \mbox{if $r=2q+1$.}
\end{cases}
\]
\end{Proposition}

\begin{Proof}
By Proposition \ref{T_4}, the second smallest reverse Wiener index
in $\mathbb{T}_{n,4}$ is precisely achieved by the smallest reverse
Wiener index of trees in $\mathbb{T}_{n,4}\setminus
\mathcal{T}_{n,4}$. Let $T$ be a tree with the smallest reverse
Wiener index in $\mathbb{T}_{n,4}\setminus \mathcal{T}_{n,4}$.

\noindent {\bf Case 1.} $T \in \mathcal{CT}_n$. Then $T$ may be
written as
$T=T_{n,k}\left(n_0;n_1^{[b_1]},n_2^{[b_2]},\ldots,n_s^{[b_s]}\right)$.
From the expression for $W(T)$ (given previous to Lemma
\ref{balace}), $\Lambda(T)$ is increasing with respect to $n_0$.
Thus $n_0=1$, $k \ge 2$ and $\sum_{i=1}^{s}{b_in_i}=n-k-2$, and by
Lemma \ref{balace}, we have $s \le 2$, and $n_2-n_1 = 1$ if $s=2$.
Using the same method in the proof of Lemma \ref{k^2+1} to analyze
$\Lambda\left(T_{n,k}\left(1;n_1^{[b_1]},n_2^{[b_2]},\ldots,n_s^{[b_s]}\right)\right)
=4n-5+(n-2)k+\sum_{i=1}^{s}{b_in_i^2}$, we have $k(T)=q-1,q$, $q+1$.
The rest proof is similar to the proof of Proposition \ref{T_4}.

If $r=1$, then $T=T_{n,q}\left(1;q-2^{[1]},q-1^{[q-1]}\right)$ with
$\Lambda(T)=2q^3+2q^2-2q+2$.

If $r=2$, then $T=T_{n,q}\left(1;q-1^{[q]}\right)$ with
$\Lambda(T)=2q^3+2q^2+q+3$.

If $r=3,\ldots,q+1$, then
$T=T_{n,q}\left(1;q-1^{[q-r+2]},q^{[r-2]}\right)$ with
$\Lambda(T)=2q^3+2q^2+3rq-5q+3r-3$.

If $r=q+2$, then $T=T_{n,q}\left(1;q^{[q]}\right)$ with
$\Lambda(T)=2q^3+5q^2+4q+3$.

If $r=q+3,\ldots,2q+1$, then
$T=T_{n,q+1}\left(1;q-1^{[2q-r+3]},q^{[r-q-2]}\right)$ with
$\Lambda(T)=2q^3+2q^2+3rq-6q+4r-4$.

\noindent {\bf Case 2.}  $T \not \in \mathcal{CT}_n$. Then $T$ may
be written as
$T=T_{n,k}\left(n_1^{[b_1]},n_2^{[b_2]},\ldots,n_s^{[b_s]}\right)$.
Note that $T\not\in {\cal T}_{n,4}$. There are two subcases.

\noindent {\bf Subcase 2.1.} $s\le 2$, and $n_2-n_1=1$ if $s=2$.

(i) $r=1$. By the monotonicity of $f(k)$ in the proof of Lemma
\ref{k^2+1} and Proposition \ref{d=3}, we have $k=q-1$, $q+1$. If
$k=q-1$, then $T=T_{n,q-1}\left(q^{[q-2]},q+1^{[1]}\right)$ with
$\Lambda(T)=2q^3+q^2+q+2$. If $k=q+1$, then
$T=T_{n,q+1}\left(q-2^{[q]},q-1^{[1]}\right)$ with
$\Lambda(T)=2q^3+q^2+q$. Thus
$T=T_{n,q+1}\left(q-2^{[q]},q-1^{[1]}\right)$.

(ii) $r=2, \dots, q$. Then $k=q-1$, $q+1$. If $k=q-1$, then
\[
T=\begin{cases}
T_{n,q-1}\left(q^{[q-r-1]},q+1^{[r]}\right)& \mbox{if } r=2,\ldots,q-2,\\
T_{n,q-1}\left(q+1^{[q-1]}\right)& \mbox{if } r=q-1,\\
T_{n,q-1}\left(q+1^{[q-2]},q+2^{[1]}\right) & \mbox{if } r=q,
\end{cases}
\]
where
\begin{eqnarray*}
& &
   \Lambda\left(T_{n,q-1}\left(q^{[q-r-1]},q+1^{[r]}\right)\right)\\
&=&
   2n(n-1)-[(n-1)(n-q)+(q+1)(n-q-1)(q-r-1)\\
   &&+(q+2)(n-q-2)r]\\
&=&
   n^2-nq^2+nq-nr+q^3+q^2+2rq-2q+3r-1\\
&=&
   2q^3+q^2+3rq-2q+3r-1,
\end{eqnarray*}
\[
\Lambda\left(T_{n,q-1}\left(q+1^{[q-1]}\right)\right)=2q^3+4q^2-2q-4,
\]
\[
\Lambda\left(T_{n,q-1}\left(q+1^{[q-2]},q+2^{[1]}\right)\right)=2q^3+4q^2+q+1.
\]
If $k=q+1$, then $T= T_{n,q+1}\left(q-2^{[q-r+1]},q-1^{[r]}\right)$,
where
\begin{eqnarray*}
& &
   \Lambda\left(T_{n,q+1}\left(q-2^{[q-r+1]},q-1^{[r]}\right)\right)\\
&=&
   2n(n-1)-\left[(n-1)(n-2-q)\right.\\
   && \left. +(q-1)(n-q+1)(q-r+1)+q(n-q)r\right]\\
&=&
   n^2+2n-nq^2+nq-nr+q^3-q^2+2rq-2q-r-1\\
&=&
   2q^3+q^2+3rq-2q+r-1.
\end{eqnarray*}
Thus $T=T_{n,q+1}\left(q-2^{[q-r+1]},q-1^{[r]}\right)$ with
$\Lambda(T)=2q^3+q^2+3rq-2q+r-1$.

(iii)  $r=q+1$. Then $k=q-1$, $q+2$. If $k=q-1$, then
$T=T_{n,q-1}\left(q+1^{[q-3]},q+2^{[2]}\right)$ with
$\Lambda(T)=2q^3+4q^2+4q+6$, and if $k=q+2$, then
$T=T_{n,q+2}\left(q-2^{[q]},q-1^{[2]}\right)$ with
$\Lambda(T)=2q^3+4q^2+4q$. Thus
$T=T_{n,q+2}\left(q-2^{[q]},q-1^{[2]}\right)$ with
$\Lambda(T)=2q^3+4q^2+4q$.

(iv) $r=q+2,\ldots,2q$. Then $k=q$, $q+2$. If $k=q$, then
$T=T_{n,q}\left(q^{[2q+1-r]},q+1^{[r-q-1]}\right)$ with
$\Lambda(T)=2q^3+q^2+3rq-5q+4r-4$ (which is shown in the proof of
Proposition \ref{T_4}), and if $k=q+2$, then
$T=T_{n,q+2}\left(q-2^{[2q+1-r]},q-1^{[r-q+1]}\right)$, where
\begin{eqnarray*}
& &
   \Lambda\left(T_{n,q+2}\left(q-2^{[2q+1-r]},q-1^{[r-q+1]}\right)\right)\\
&=&
   2n(n-1)-[(n-1)(n-3-q)\\
   && +(q-1)(n-q+1)(2q+1-r)+q(n-q)(r-q-1)]\\
&=&
   2q^3+q^2+3rq-q+2r-2\\
&>&\Lambda\left(T_{n,q}\left(q^{[2q+1-r]},q+1^{[r-q-1]}\right)\right).
\end{eqnarray*}
Thus $T=T_{n,q}\left(q^{[2q+1-r]},q+1^{[r-q-1]}\right)$ with
$\Lambda(T)=2q^3+q^2+3rq-5q+4r-4$.

(v) $r=2q+1$. Then $k=q$, $q+2$. If $k=q$, then
$T=T_{n,q}\left(q+1^{[q]}\right)$, and if $k=q+2$, then
$T=T_{n,q+2}\left(q-1^{[q+2]}\right)$, both with reverse Wiener
index $2q^3+7q^2+6q$.

\noindent {\bf Subcase 2.2.} There exist $1 \le i<j \le s$ such that
$n_j-n_i \ge 2$. Let $T^*$ be the tree constructed as in the proof
of Lemma \ref{balace} from $T$. By the proof there, we have
$\Lambda(T) = \Lambda(T^*)+2(n_j-n_i-1)$. Thus $\Lambda(T)$ is
minimum if and only if $T^*\in \mathcal{T}_{n,4}$ and $n_j-n_i=2$.

(i) $r=1$. Then
$T=T_{n,q}\left(q-2^{[1]},q-1^{[q-2]},q^{[r]}\right)$ with
$\Lambda(T)=2q^3+q^2+2$.

(ii) $r=2,\ldots,q$. Then $\Lambda(T)$ is at least the minimum value
of the reverse Wiener indices of
\begin{eqnarray*}
T_{n,q}\left(q-2^{[1]},q-1^{[q-3]},q^{[2]}\right),
T_{n,q}\left(q-2^{[1]},q-1^{[q-2]},q+1^{[1]}\right) \mbox{ if }
r=2,\\
 T_{n,q}\left(q-1^{[q-1]},q+1^{[1]}\right),
 T_{n,q}\left(q-2^{[1]},q-1^{[q-4]},q^{[3]}\right), \\
 T_{n,q}\left(q-2^{[1]},q-1^{[q-3]},q^{[1]},q+1^{[1]}\right)
 \mbox{ if } r=3,\\
T_{n,q}\left(q-1^{[q-r+2]},q^{[r-3]},q+1^{[1]}\right),
T_{n,q}\left(q-2^{[1]},q-1^{[q-r-1]},q^{[r]}\right), \\
T_{n,q}\left(q-2^{[1]},q-1^{[q-r]},q^{[r-2]},q+1^{[1]}\right)
 \mbox{ if } r=4,\ldots,q-2,\\
T_{n,q}\left(q-1^{[3]},q^{[q-4]},q+1^{[1]}\right),
T_{n,q}\left(q-2^{[1]},q^{[q-1]}\right), \\
T_{n,q}\left(q-2^{[1]},q-1^{[1]},q^{[q-3]},q+1^{[1]}\right)
 \mbox{ if } r=q-1,\\
T_{n,q}\left(q-1^{[2]},q^{[q-3]},q+1^{[1]}\right),
T_{n,q}\left(q-2^{[1]},q^{[q-2]},q+1^{[1]}\right)
 \mbox{ if } r=q.
\end{eqnarray*}
By comparing the reverse Wiener indices of these trees, we have
\[
T=
\begin{cases}
T_{n,q}\left(q-2^{[1]},q-1^{[q-r-1]},q^{[r]}\right) & \mbox{if }
r=2,\ldots,q-2,\\
T_{n,q}\left(q-2^{[1]},q^{[q-1]}\right) & \mbox{if }
r=q-1,\\
T_{n,q}\left(q-2^{[1]},q^{[q-2]},q+1^{[1]}\right) & \mbox{if } r=q,
\end{cases}
 \] where
\begin{eqnarray*}
& &
   \Lambda\left(T_{n,q}\left(q-2^{[1]},q-1^{[q-r-1]},q^{[r]}\right)\right)\\
&=&
   2n(n-1)-[(n-1)(n-q-1)+(q-1)(n-q+1)\\
&&+q(n-q)(q-r-1)+(q+1)(n-q-1)r]\\
&=&
   n^2+n-nq^2+nq-nr+q^3+2rq-3q+r\\
&=&
   2q^3+q^2+3rq-3q+2r,
\end{eqnarray*}
\[
\Lambda\left(T_{n,q}\left(q-2^{[1]},q^{[q-1]}\right)\right)=2q^3+4q^2-4q-2,
\]
\[
\Lambda\left(T_{n,q}\left(q-2^{[1]},q^{[q-2]},q+1^{[1]}\right)\right)=2q^3+4q^2-q+2.
\]

(iii) $r=q+1$. Then
$T=T_{n,q}\left(q-1^{[1]},q^{[q-2]},q+1^{[1]}\right)$,
$T_{n,q+1}\left(q-2^{[1]},q-1^{[q-1]},q^{[1]}\right)$, both with
reverse Wiener index $2q^3+4q^2+2q+2$.

(iv) $r=q+2,\ldots,2q$. Then $\Lambda(T)$ is at least the minimum
value of the reverse Wiener indices of
\begin{eqnarray*}
T_{n,q+1}\left(q-2^{[1]},q-1^{[q-2]},q^{[2]}\right),
T_{n,q+1}\left(q-2^{[1]},q-1^{[q-1]},q+1^{[1]}\right) \mbox{ if }
r=q+2,\\
T_{n,q+1}\left(q-2^{[1]},q-1^{[q-3]},q^{[3]}\right),
T_{n,q+1}\left(q-2^{[1]},q-1^{[q-2]},q^{[1]},q+1^{[1]}\right),\\
T_{n,q+1}\left(q-1^{[q]},q+1^{[1]}\right) \mbox{if } r=q+3,\\
T_{n,q+1}\left(q-2^{[1]},q-1^{[2q-r]},q^{[r-q]}\right),
T_{n,q+1}\left(q-2^{[1]},q-1^{[2q+1-r]},q^{[r-q-2]},q+1^{[1]}\right),\\
T_{n,q+1}\left(q-1^{[2q+3-r]},q^{[r-q-3]},q+1^{[1]}\right) \mbox{if
} r=q+4,\ldots,2q-1,\\
T_{n,q+1}\left(q-2^{[1]},q^{[q]}\right),
T_{n,q+1}\left(q-2^{[1]},q-1^{[1]},q^{[q-2]},q+1^{[1]}\right),\\
T_{n,q+1}\left(q-1^{[3]},q^{[q-3]},q+1^{[1]}\right) \mbox{if } r=2q.
\end{eqnarray*}
Thus we have
\[
T=
\begin{cases}
T_{n,q+1}\left(q-2^{[1]},q-1^{[q-2]},q^{[2]}\right) \mbox{~~~~~~}
\mbox{if $r=q+2$,}\\

T_{n,q+1}\left(q-2^{[1]},q-1^{[q-3]},q^{[3]}\right),
T_{n,q+1}\left(q-1^{[q]},q+1^{[1]}\right) \\
\mbox{~~~~~~~~~~~~~~~~~~~~~~~~~~~~~~~~~~~~~~~~~~~~~}  \mbox{if $r=q+3$,}\\

T_{n,q+1}\left(q-2^{[1]},q-1^{[2q-r]},q^{[r-q]}\right),
T_{n,q+1}\left(q-1^{[2q+3-r]},q^{[r-q-3]},q+1^{[1]}\right) \\
\mbox{~~~~~~~~~~~~~~~~~~~~~~~~~~~~~~~~~~~~~~~~~~~~~}
\mbox{if $r=q+4,\ldots,2q-1$,}\\

T_{n,q+1}\left(q-2^{[1]},q^{[q]}\right),
T_{n,q+1}\left(q-1^{[3]},q^{[q-3]},q+1^{[1]}\right) \\
\mbox{~~~~~~~~~~~~~~~~~~~~~~~~~~~~~~~~~~~~~~~~~~~~~} \mbox{if
$r=2q$,}
\end{cases}
\]
where
\[
\Lambda\left(T_{n,q+1}\left(q-2^{[1]},q-1^{[q-2]},q^{[2]}\right)\right)=2q^3+4q^2+5q+5,
\]
\begin{eqnarray*}
   \Lambda\left(T_{n,q+1}\left(q-2^{[1]},q-1^{[q-3]},q^{[3]}\right)\right)
&=&
   \Lambda\left(T_{n,q+1}\left(q-1^{[q]},q+1^{[1]}\right)\right)\\
&=&
   2q^3+4q^2+8q+8,
\end{eqnarray*}
\begin{eqnarray*}
& &
   \Lambda\left(T_{n,q+1}\left(q-2^{[1]},q-1^{[2q-r]},q^{[r-q]}\right)\right)\\
&=&
   \Lambda\left(T_{n,q+1}\left(q-1^{[2q+3-r]},q^{[r-q-3]},q+1^{[1]}\right)\right)\\
&=&
   2n(n-1)-[(n-1)(n-2-q)+(q-1)(n-q+1)+q(n-q)(2q-r)\\
   &&+(q+1)(n-q-1)(r-q)]\\
&=&
   n^2+2n-nq^2+nq-nr+q^3-q^2+2rq-4q+r-1\\
&=&
   2q^3+q^2+3rq-4q+3r-1,
\end{eqnarray*}
\[
\Lambda\left(T_{n,q+1}\left(q-2^{[1]},q^{[q]}\right)\right)=
\Lambda\left(T_{n,q+1}\left(q-1^{[3]},q^{[q+3]},q+1^{[1]}\right)\right)=2q^3+7q^2+2q-1.
\]

(v) $r=2q+1$. Then $\Lambda(T)$ is at least the minimum value of the reverse Wiener indices of
$T_{n,q+1}\left(q-2^{[1]}, q^{[q-1]},q+1^{[1]}\right)$
and
$T_{n,q+1}\left(q-1^{[2]}, q^{[q-2]},q+1^{[1]}\right)$.
By direct calculation, we have $T=
T_{n,q+1}\left(q-1^{[2]},q^{[q-2]},q+1^{[1]}\right)$ with
$\Lambda(T)=2q^3+7q^2+5q+2$.

By comparing the above cases, we have: if $r=1,\ldots,q-1$, then
$\Lambda(T)$ is precisely the minimum value of
$2q^3+2q^2+3rq-5q+3r-3$, $2q^3+q^2+3rq-2q+r-1$ and
$2q^3+q^2+3rq-3q+2r$; if $r=q$, then $\Lambda(T)$ is precisely the
minimum value of $2q^3+5q^2-2q-3$, $2q^3+4q^2-q-1$ and
$2q^3+4q^2-q+2$; if $r=q+1$, then $\Lambda(T)$ is precisely the
minimum value of $2q^3+5q^2+q$, $2q^3+4q^2+4q$ and $2q^3+4q^2+2q+2$;
if $r=q+2$, then $\Lambda(T)$ is precisely the minimum value of
$2q^3+5q^2+4q+3$, $2q^3+4q^2+5q+4$ and $2q^3+4q^2+5q+5$; if
$r=q+3,\ldots,2q+1$, then $\Lambda(T)$ is precisely the minimum
value of $2q^3+2q^2+3rq-6q+4r-4$, $2q^3+q^2+3rq-5q+4r-4$ and
$2q^3+q^2+3rq-4q+3r-1$. Now the result follows easily.
\end{Proof}

\begin{Theorem}
Among the trees in $\mathcal{T}_n$ with $n \ge 5$, the trees in
$\mathcal{T}'_{n,4}$ for $n=5$ or $n \ge 58$,
$D_{n,\lfloor\frac{n}{2}\rfloor-1}$ for $6 \le n \le 57$ are the
unique trees with the third smallest reverse Wiener index, which is
equal to $20$ for $n=5$,
$\frac{n^2}{2}+\frac{3n}{2}-2-\left(\lfloor\frac{n}{2}\rfloor-1\right)
\left(\lceil\frac{n}{2}\rceil+1\right)$ for $6\le n\le 57$, and
$g(n,4)$ for $n\ge 58$, where $\mathcal{T}'_{n,4}$ and $g(n,4)$ are
given in Proposition \ref{T_4-sec}.
\end{Theorem}

\begin{Proof}
The case $n=5$ is trivial. Suppose that $n\ge 6$. By Theorem
\ref{sec-min}, among the trees in $\mathcal{T}_n$ with $n \ge 6$,
$D_{n,\lfloor\frac{n}{2}\rfloor}$ for $6 \le n \le 56$, $D_{57,28}$,
$T_{57,7}\left(7^{[7]}\right)$ and $T_{57,8}\left(6^{[8]}\right)$
for $n=57$, and the trees in $\mathcal{T}_{n,4}$ for $n \ge 58$ are
the unique trees with the second smallest reverse Wiener index,
which are equal to $f(n,3)$ for $n\le 56$, $f(n,3)=f(n,4)$ for
$n=57$ and $f(n,4)$ for $n\ge 58$. By Proposition \ref{T_4-sec}, and
Lemmas \ref{d-2} and \ref{T-5}, we have $f(n,d)
> g(n,4)$ for $5 \le d \le n-1$. Thus the third smallest reverse
Wiener index in $\mathcal{T}_n$ is equal to $\min\{g(n,3), f(n,4)\}$
for $6 \le n \le 56$, $\min\{g(n,3),g(n,4)\}$ for $n=57$ and
$\min\{f(n,3), g(n,4)\}$ for $n \ge 58$, and by Propositions
\ref{d=3} and \ref{T_4-sec}, it is precisely achieved by graphs in
$\mathcal{T'}_{n,3}$ and $\mathcal{T}_{n,4}$ for $6 \le n \le 56$,
$\mathcal{T'}_{n,3}$ and $\mathcal{T'}_{n,4}$ for $n=57$, and
$\mathcal{T}_{n,3}$ and $\mathcal{T'}_{n,4}$ for $n\ge 58$.

Recall that the expressions for $f(n,3)$ and $g(n,3)$ are given in
Proposition \ref{d=3}, while the expressions for $f(n,4)$ and
$g(n,4)$ are given in Propositions \ref{T_4} and \ref{T_4-sec},
respectively. Let $n$ be of the form $(\ref{formu-n})$.

 By direct checking, we have $g(n, 3)<f(n,4)$  for $6 \le
n \le 56$, and $g(57,3)<g(57,4)$. Thus the results for cases
$n=6,\dots,57$ follow from Proposition \ref{d=3}. Suppose that $n
\ge 58$.

\noindent {\bf Case 1.} $r=1,\ldots,q-1$. Since $q^2+1 \ge 58$, we
have $q \ge 8$. Note that $\frac{1}{2}q^2-3q-\frac{1}{2}>0$. Then
\begin{eqnarray*}
& &
   f(n,3)-g(n,4)\\
&\ge&
   \frac{1}{4}q^4+\frac{1}{2}rq^2+\frac{3}{2}q^2+\frac{1}{4}r^2+\frac{3}{2}r-2
   -\left(2q^3+q^2+3rq-3q+2r\right)\\
&=&
   \frac{1}{4}q^4-2q^3+\frac{1}{2}q^2+3q-2+\frac{1}{4}r^2+\left(\frac{1}{2}q^2-3q-\frac{1}{2}\right)r\\
&\ge&
   \frac{1}{4}q^4-2q^3+\frac{1}{2}q^2+3q-2+\frac{1}{4}+\frac{1}{2}q^2-3q-\frac{1}{2}\\
&=&
   \frac{1}{4}q^4-2q^3+q^2-\frac{9}{4}>0,
\end{eqnarray*}
and thus $f(n,3)>g(n,4)$ for $n \ge 65$.

%
%
%
%
%
%
%
%
%

\noindent {\bf Case 2.}  $r=q+2$. Since $q^2+q+2 \ge 58$, we have $q
\ge 7$. Then
\[
f(n,3)-g(n,4) \ge
\frac{1}{4}q^4-\frac{3}{2}q^3-\frac{5}{4}q^2-\frac{5}{2}q-2>0,
\]
and thus $f(n,3)>g(n,4)$ for $n \ge 58$.

\noindent {\bf Case 3.} $r=q, q+1,\ldots,2q+1$.
Note that $\frac{1}{2}q^2-3q-\frac{3}{2}>0$. We
have
\begin{eqnarray*}
& &
   f(n,3)-g(n,4)\\
&\ge&
   \frac{1}{4}q^4+\frac{1}{2}rq^2+\frac{3}{2}q^2+\frac{1}{4}r^2+\frac{3}{2}r-2
   -\left(2q^3+q^2+3rq-4q+3r-1\right)\\
&=&
   \frac{1}{4}q^4-2q^3+\frac{1}{2}q^2+4q-1+\frac{1}{4}r^2+\left(\frac{1}{2}q^2-3q-\frac{3}{2}\right)r.
\end{eqnarray*}
If $r=q, q+1$, then $q\ge 8$, and thus $f(n,3)-g(n,4)\ge \frac{1}{4}q^4-\frac{3}{2}q^3-\frac{9}{4}q^2+\frac{5}{2}q-1>0$, while if $r=q+3,...,2q+1$, then $q\ge 7$, and thus $f(n,3)-g(n,4)\ge \frac{1}{4}q^4-\frac{3}{2}q^3-\frac{3}{4}q^2-5q-\frac{13}{4}>0$.
Thus $f(n,3)>g(n,4)$ for $n \ge 59$.

Combining Cases 1--3 above, we have $f(n,3) > g(n,4)$ for $n \ge
58$. Then the result follows from Proposition \ref{d=3}.
\end{Proof}

\vspace{4mm}

\noindent {\it Acknowledgement.\/} This work was supported by the
National Natural Science Foundation of China (no.~11071089).


\begin{thebibliography}{99}

\bibitem{BaMiIB} A. T. Balaban, D. Mills, O. Ivanciuc, S. C. Basak,
Reverse Wiener indices, Croat. Chem. Acta 73 (2000)
923--941.


\bibitem{CaZh}
X. Cai, B. Zhou, Reverse Wiener indices of connected graphs,
 MATCH Commun. Math. Comput. Chem. 60 (2008) 95--105.


\bibitem{DoEnG}
A. A. Dobrynin, R. Entringer,  I. Gutman,  Wiener index of trees:
theory and applications, Acta Appl. Math.  66 (2001)
211--249.

\bibitem{GKM1} I. Gutman, S. Klav\v zar, B. Mohar (Eds.), Fifty years of the Wiener index,
MATCH Commun. Math. Comput. Chem. 35 (1997) 1--259.


\bibitem{GKM2}I. Gutman, S. Klav\v zar, B. Mohar (Eds.), Fiftieth anniversary of the
Wiener index, Discrete Appl. Math. 80 (1997) 1--113.

\bibitem{GP} I. Gutman, O. E. Polansky,
Mathematical Concepts in Organic Chemistry,
Springer--Verlag, Berlin 1986.


\bibitem{Ho} H. Hosoya, Topological index.
A newly proposed quantity characterizing the topological nature of
structural isomers of saturated hydrocarbons, Bull. Chem. Soc. Japan 44 (1971) 2332--2339.

\bibitem{IIB}  O. Ivanciuc,  T. Ivanciuc, A. T. Balaban,
Quantitative structure--property relationship evaluation of
structural descriptors derived from the distance and reverse
Wiener matrices, Internet Electron. J. Mol. Des. 1
(2002) 467--487.


\bibitem{Kar}
M. Karelson, Molecular Descriptors in QSAR/QSPR, Wiley,
New York, 2000.

\bibitem{LZ1}
W. Luo, B. Zhou, Further properties of reverse Wiener index,
MATCH Commun. Math. Comput. Chem. 61 (2009) 653--661.

\bibitem{LZ2} W. Luo, B. Zhou, On ordinary and reverse Wiener indices of non-caterpillars,
Math. Comput. Modelling 50 (2009) 188--193.

\bibitem{Pl}J. Plesn\'ik, On the sum of all distances in a graph or digraph, J. Graph
Theory 8 (1984) 1--21.


\bibitem{Rou} D. H. Rouvray,  The rich legacy of half a century of the Wiener index, in:   D. H.
Rouvray, R. B. King (Eds.),Topology
in Chemistry -- Discrete Mathematics of Molecules, Horwood, Chichester 2002, pp. 16--37.

\bibitem{ToCo} R. Todeschini, V. Consonni,
Handbook of Molecular Descriptors, Wiley--VCH, Weinheim
2000.

\bibitem{Tr} N. Trinajsti\'c, Chemical Graph Theory, 2nd revised ed., CRC press, Boca Raton, 1992.


\bibitem{Wie} H. Wiener, Structural determination of paraffin boiling points, J. Am.
Chem. Soc. 69 (1947) 17--20.

\bibitem{ZhZh} B. Zhang, B. Zhou, Modified and reverse Wiener
indices of trees, Z. Naturforsch. 61a (2006)
536--540.

\bibitem{Zh}
B. Zhou, Reverse Wiener index, in: I. Gutman, B. Furtula (Eds.), Novel Molecular Structure Descriptors -- Theory and Applications II, Univ. Kragujevac, Kragujevac, 2010, pp. 193--204.

\end{thebibliography}
\end{document}